\documentclass[11pt]{amsart}
\usepackage{amsmath,amssymb,euscript,mathrsfs}
\usepackage[bookmarks=false]{hyperref}

\newtheorem{thm}{Theorem}
\newtheorem{lem}[thm]{Lemma}
\newtheorem{propn}[thm]{Proposition}

\theoremstyle{definition}
\newtheorem{rmk}[thm]{Remark}

\newcommand{\isom}{\cong}
\DeclareMathOperator{\rk}{rk}
\DeclareMathOperator{\Aut}{Aut}

\renewcommand{\P}{\mathbb{P}}
\newcommand{\Z}{\mathbb{Z}}
\newcommand{\C}{\mathbb{C}}

\renewcommand{\O}{\mathcal{O}}

\renewcommand{\phi}{\varphi}
\renewcommand{\setminus}{\smallsetminus}

\renewcommand{\tilde}{\widetilde}

\newcommand{\Fl}{Fl}
\newcommand{\FFl}{\mathbf{Fl}}
\newcommand{\OOmega}{\mathbf{\Omega}}

\begin{document}

\title[Positivity in the cohomology of flag bundles]{Positivity in the cohomology of flag bundles (after Graham)}
\author{Dave Anderson}
\date{September 30, 2007}
\maketitle

In \cite{graham}, Graham proves that the structure constants of the equivariant cohomology ring of a flag variety are positive combinations of monomials in the roots:
\begin{thm}[{\cite[Cor.\ 4.1]{graham}}]
Let $X = G/B$ be the flag variety for a complex semisimple group $G$ with maximal torus $T\subset B$, and let $\{ \sigma_w \in H_T^*X \,|\, w\in W\}$ be the basis of ($B$-invariant) Schubert classes.  Let $\{\alpha_i\}$ be the simple roots which are \emph{negative} on $B$.  Then in the expansion
\begin{eqnarray*}
\sigma_u\cdot\sigma_v = \sum_w c_{uv}^w \, \sigma_w,
\end{eqnarray*}
the coefficients $c_{uv}^w$ are in $\Z_{\geq 0}[\alpha]$.
\end{thm}

\noindent
Graham deduces this from a more general result about varieties with finitely many unipotent orbits, which is proved using induction and a calculation in the rank-one case.

The goal of this note is to give a short, geometric proof of Graham's positivity theorem, based on a transversality argument.  Here I only discuss type $A$, but other types work as well.  (For a type-uniform version, a change of language is needed: one should replace vector bundles with corresponding principal $G$-bundles.)

Throughout, $\Fl$ denotes the variety of (complete) flags in $\C^n$, and if $V \to X$ is a vector bundle, $\FFl(V) \to X$ is the bundle of flags in $V$.

Recall that for $T' \isom (\C^*)^{n}$, we have $BT' = (\P^\infty)^{\times n}$ and $H_{T'}^*\Fl = H^*(ET'\times^{T'} \Fl) = H^*\FFl(E')$, where $E'$ is the sum of the $n$ tautological line bundles on $BT'$.  The \emph{effective} action on $\Fl$ is by $T \isom (\C^*)^n/\C^*$, and the classifying space for this torus is $BT=(\P^\infty)^{\times n-1}$.  We will usually deal with the effective torus.

Let $\P = \P^m \times \cdots  \times \P^m$ ($n-1$ factors), with $m\gg 0$, and write $H^*\P = \Z{[\alpha_1,\ldots,\alpha_{n-1}]}$.  (We always assume that $m$ is large enough so that there are no relations in the relevant degrees.)  Let $M_i = p_1^*(\O(-1))$ be the tautological bundle on the $i$th factor, and let $\alpha_i = -c_1(M_i)$.  Note that the class of any effective cycle in $H^*\P$ is a positive polynomial in the $\alpha$'s.

Let 
\begin{eqnarray*}
L_i = M_1 \otimes \cdots \otimes M_{i-1}
\end{eqnarray*}
for $1\leq i\leq n$ (so $L_1 = \O$ is the trivial line bundle), and let $E_i = L_1 \oplus \cdots \oplus L_i$.  Thus we have a flag $E_{\bullet}$ in $E = E_n$.  Let $\tilde{E}_{\bullet}$ be the opposite flag, with $\tilde{E}_i = L_n \oplus \cdots \oplus L_{n+1-i}$.  In the flag bundle $p:\FFl(E) \to \P$, with universal quotient flags $Q_{\bullet}$, we have Schubert loci $\OOmega_w = \Omega_w(E_{\bullet} \to Q_{\bullet})$, defined by
\begin{eqnarray}\label{degloci}
\OOmega_w = \{ x\in \FFl(E) \,|\, \rk(E_p \to Q_q) \leq \#(i\leq q\,|\, w(i)\leq p) \}.
\end{eqnarray}
Opposite Schubert loci $\tilde\OOmega_w = \Omega_w(\tilde{E}_{\bullet} \to Q_\bullet)$ are defined similarly.  We also have ``Schubert cell bundles'' $\OOmega^o_w$: these are affine bundles over $\P$ which are open in the corresponding loci $\OOmega_w$, and are defined by replacing the inequality in \eqref{degloci} with an equality.

The classes ${[\OOmega_w]}$ form a basis for $H^*\FFl(E)$ over $H^*\P$, as $w$ ranges over $S_n$.  Writing 
\begin{eqnarray*}
{[\OOmega_u]}\cdot{[\OOmega_v]} = \sum_w c_{uv}^w {[\OOmega_w]}
\end{eqnarray*}
with $c_{uv}^w \in H^*\P$, our main result is the following:

\begin{propn} \label{prop:positive}
The polynomials $c_{uv}^w$ are \emph{positive}, that is, $c_{uv}^w \in \Z_{\geq 0}{[\alpha_1,\ldots,\alpha_{n-1}]}$.
\end{propn}

This implies Graham's positivity theorem (in this context), since $\P$ approximates $BT$ for $m$ sufficiently large, and $\FFl(E)$ approximates $ET\times^T \Fl$, with $[\OOmega_w]$ corresponding to the equivariant class $\sigma_w$.  (See \cite[\S9]{eq}.)

Proposition \ref{prop:positive} is a consequence of a transversality statement:

\begin{propn} \label{prop:transverse}
For any $u,v,w \in S_n$, there is a translate $\OOmega'_v$ of $\OOmega_v$ by the action of a connected algebraic group such that $\OOmega'_v$ intersects $\OOmega_u$ and $\tilde\OOmega_{w_0\,w}$ properly and generically transversally.
\end{propn}

To deduce Proposition \ref{prop:positive}, first note that the intersection $\OOmega_u \cap \tilde\OOmega_{w_0\,w}$ is always proper and generically transverse.  Thus Proposition \ref{prop:transverse} says that $\OOmega'_v \cap (\OOmega_u \cap \tilde\OOmega_{w_0\,w})$ is proper and generically transverse.    By \cite[Ex.\ (8.1.11)]{it},
 this says that 
\begin{eqnarray*}
[\OOmega_v]\cdot [\OOmega_u] \cdot [\tilde\OOmega_{w_0\,w}] = [\OOmega'_v \cap \OOmega_u \cap \tilde\OOmega_{w_0\,w}].
\end{eqnarray*}
(Since $\OOmega'_v = g\cdot \OOmega_v$ for some $g$ in a connected algebraic group, $[\OOmega'_v] = [\OOmega_v]$.)  Using relative Poincar\'e duality (see e.g. \cite[\S A.6]{eq}), we have
\begin{eqnarray*}
c_{uv}^w = p_*([\OOmega_u]\cdot[\OOmega_v]\cdot[\tilde\OOmega_{w_0 w}]) = p_*([\OOmega_u \cap \OOmega'_v \cap \tilde\OOmega_{w_0\,w}]).
\end{eqnarray*}
This is an effective class in $H^*\P$, so Proposition \ref{prop:positive} follows.

\begin{proof}[Proof of Proposition \ref{prop:transverse}]
This is essentially an application of Kleiman's theorem.  The endomorphism bundle
\begin{eqnarray*}
\mathbf{End}(E) &=& \bigoplus_{i,j} L_i^{-1} \otimes L_j \\
&=& \left(\bigoplus_{i<j} M_i\otimes\cdots\otimes M_{j-1}\right) \oplus \O^{\oplus n} \oplus \left(\bigoplus_{i>j} M_j^{-1} \otimes \cdots \otimes M_{i-1}^{-1} \right)
\end{eqnarray*}
has global sections in lower-triangular matrices, so the group $B$ of (invertible) lower-triangular matrices acts on $\FFl(E)$, fixing the flag $\tilde{E}_{\bullet}$ and stabilizing $\tilde\OOmega_{w_0\, w}$.  (Note that the entries of a matrix in $B$ are global sections of the line bundles $M_j^{-1} \otimes \cdots \otimes M_{i-1}^{-1}$, i.e., multi-homogeneous polynomials.  This is a connected group over $\C$, acting on a fiber $p^{-1}(x) \subset \FFl(E)$ by first evaluating the sections at $x$.)

Now let $H = (GL_{m+1})^{\times (n-1)}$, and for $b\in B$, let $b_x$ be the evaluation at $x\in \P$ (so the action of $b$ on $p^{-1}(x)$ is by $b_x$). Consider the semidirect product $\Gamma = B\rtimes H$, given by $(h\cdot b \cdot h^{-1})_x = b_{h^{-1}\cdot x}$.  (This action of $H$ on $B$ is just the usual action of $H$ on global sections of the equivariant vector bundle $\mathbf{End}(E)$.)\footnote{Alternatively, one could take $\Gamma$ to be the subgroup of $\Aut(\FFl(E))$ generated by the images of $B$ and $G$ via the homomorphisms corresponding to their respective actions.}  As a semidirect product of connected groups, $\Gamma$ is a connected algebraic group.  
We claim that the locus $\tilde\OOmega^o_{w_0\,w}$ is homogeneous for the action of $\Gamma$.  Indeed, $B$ acts transitively on each fiber of $\tilde\OOmega^o_{w_0\,w}$, and the action of $H$ on $\FFl(E)$ induces a transitive action on the set of fibers of $\tilde\OOmega^o_{w_0\,w}$.  (The line bundles $L_i$ are equivariant for $H$, so $H$ preserves the flag $\tilde{E}_\bullet$, and therefore acts on $\tilde\OOmega_{w_0\,w}$.)

Finally, note that $\OOmega^o_u$ and $\tilde\OOmega^o_{w_0\,w}$ intersect transversally, as do $\OOmega^o_v$ and $\tilde\OOmega^o_{w_0\,w}$.  The proposition follows from Lemma \ref{lem:kleiman} below, taking $U = \OOmega_u$, $V = \OOmega_v$, and $W = \tilde\OOmega_{w_0\,w}$, with their stratifications by Schubert loci.
\end{proof}

\begin{lem} \label{lem:kleiman}
Let $X$ be a nonsingular variety over a field of characteristic $0$, with an action of a connected algebraic group $\Gamma$.  Let $U,V,W\subset X$ be subvarieties with stratifications
\begin{eqnarray*}
U_0 \subset \cdots \subset U_\ell = U, \\
V_0 \subset \cdots \subset V_m = V, \\
W_0 \subset \cdots \subset W_n = W,
\end{eqnarray*}
with each stratum $U_i \setminus U_{i-1}$ nonsingular.  Assume also that $\Gamma$ acts on $W$, with each stratum $W_i\setminus W_{i-1}$ a disjoint union of homogeneous spaces.

If $U_i\setminus U_{i-1}$ meets $W_k\setminus W_{k-1}$ transversally for all $i,k$, and similarly for $V_j\setminus V_{j-1}$ and $W_k\setminus W_{k-1}$, then there is an element $g\in \Gamma$ such that $g\cdot V$ meets $U\cap W$ properly and generically transversally.
\end{lem}

\noindent
This can be deduced from results found in \cite{speiser}; see also \cite{sierra} for a vast generalization.  The proof of this version is quite short, so we give it here.

\begin{proof}
Applying Kleiman's theorem (cf.\ \cite[III.10.8]{hartshorne}) to the pairs $(U_i\setminus U_{i-1} \cap W_k\setminus W_{k-1})$ and $(V_j\setminus V_{j-1} \cap W_k\setminus W_{k-1})$ inside the homogeneous space $W_k\setminus W_{k-1}$, we can choose $g\in \Gamma$ such that each intersection 
\begin{eqnarray*}
& &(U_i\setminus U_{i-1} \cap W_k\setminus W_{k-1}) \cap g\cdot(V_j\setminus  V_{j-1} \cap W_k\setminus W_{k-1}) \\
& &=
(U_i\setminus U_{i-1} \cap W_k\setminus W_{k-1}) \cap (g\cdot V_j\setminus g\cdot V_{j-1} \cap W_k\setminus W_{k-1})
\end{eqnarray*}
is transverse, so the intersection $U\cap W \cap g\cdot V$ is proper and generically transverse.
\end{proof}

\begin{rmk}
All that is required in the proof of Proposition \ref{prop:transverse} are the facts that $\P$ is homogeneous for the action of an algebraic group $H$, and $L_i$ are $H$-equivariant line bundles such that $L_i^{-1}\otimes L_j$ is globally generated for $i>j$.
\end{rmk}

\begin{rmk}
To recover the result that for (type $A$) equivariant Schubert calculus, the structure constants $c_{uv}^w$ are in $\Z_{\geq 0}{[t_2-t_1,\ldots,t_n-t_{n-1}]}$, let $\P' = (\P^{m'})^{\times n}$ and choose a map
$\phi: \P' \to \P$
such that $\phi^*M_i = M'_i \otimes (M'_{i+1})^{-1}$, where $M'_i$ is the tautological bundle on the $i$th factor of $\P'$, with $t_i = c_1(M'_i)$.  (Note that $\phi$ will not be holomorphic!)

The $T'$-equivariant class of a Schubert variety (for $T'=(\C^*)^n$) can be identified with the class of the locus $\Omega_w(E'_\bullet \to Q_{\bullet}) \subset \FFl(E')$, where $E'_i = M'_1 \oplus \cdots \oplus M'_i$ is a flag of bundles on $\P'$.  Since this is $\phi^{-1}\OOmega_w$, the equivariant structure constants are $\phi^*c_{uv}^w$, which are positive in the variables $\phi^*\alpha_i = t_{i+1}-t_i$.
\end{rmk}

\begin{rmk}
The naive choice of flag, with $F_i = M_1 \oplus \cdots \oplus M_i$, does not work: The bundle $\mathbf{End}(F)$ has only diagonal global sections, so the corresponding loci $\OOmega_w^o$ are not homogeneous.  This explains why one does not see positivity over $\P'$.
\end{rmk}

\noindent
{\em Acknowledgements.}  This proof was inspired by William Fulton's lectures on equivariant cohomology \cite{eq}, and I thank him for comments on the manuscript.  Thanks also to Sue Sierra for interesting discussions, and for bringing \cite{speiser} to my attention.

\small{\textsc{Department of Mathematics, University of Michigan, Ann Arbor, MI 
48109}

\textit{E-mail address}: \texttt{dandersn@umich.edu}}

\end{document}